\documentclass[12pt,a4paper]{amsart}

\usepackage{amsmath}
\usepackage{amscd}
\numberwithin{equation}{section}

\input{prepictex}
\input{pictex}
\input{postpictex}
\chardef\bslash=`\\ 

\makeatletter
\def\verbatim{\interlinepenalty\@M \@verbatim
  \leftskip\@totalleftmargin\advance\leftskip2pc
  \frenchspacing\@vobeyspaces \@xverbatim}
\makeatother
\newtheorem{theorem}{Theorem}[section]
\newtheorem{corollary}[theorem]{Corollary}
\newtheorem{lemma}[theorem]{Lemma}
\newtheorem{proposition}[theorem]{Proposition}

\newtheorem{example}[theorem]{Example}

\theoremstyle{definition}
\newtheorem{definition}[theorem]{Definition}
\newtheorem{remark}[theorem]{Remark}

\newcounter{picture}

\newsymbol\rtimes 226F 

\DeclareMathOperator{\supp}{supp}
\DeclareMathOperator{\Spec}{Sp}

\newcommand{\FF}{{\mathbb F}}

\newcommand{\PP}{{\mathbb P}}
\newcommand{\QQ}{{\mathbb Q}}
\newcommand{\RR}{{\mathbb R}}

\newcommand{\ZZ}{{\mathbb Z}}
\newcommand{\cA}{{\mathcal A}}

\newcommand{\cO}{{\mathcal O}}

\newcommand{\D}{{\Delta}}

\newcommand{\G}{{\Gamma}}
\newcommand{\Om}{{\Omega}}
\newcommand{\e}{{\epsilon}}

\newcommand{\w}{{\omega}}
\newcommand{\wa}{{\varpi}}

\newcommand{\Ind}{{\bf 1}} 
\newcommand{\PGL}{{\text{\rm{PGL}}}}
\newcommand{\PSL}{{\text{\rm{PSL}}}}

\newcommand{\tA}{{\widetilde A}}

\begin{document}

\title[$n$-filling actions]
{Simple purely infinite $C^*$-algebras and $n$-filling actions}

\date{January 14, 2000}

\author[Paul Jolissaint]{Paul Jolissaint}
\address{Institut de Math\'ematiques, Universit\'e de Neuch\^{a}tel,
Emile-Argand 11, CH-2000 Neuch\^{a}tel, Switzerland}
\email{paul.jolissaint@maths.unine.ch}

\author[Guyan Robertson]{Guyan Robertson}
\address{Department of Mathematics  \\
        University of Newcastle\\  NSW  2308\\ AUSTRALIA}
\email{guyan@maths.newcastle.edu.au }


\thanks{This research was supported by a research visitor grant from the University of Newcastle.} 
\thanks{ \hfill Typeset by  \AmS-\LaTeX}
\begin{abstract}
Let $n$ be a positive integer.
We introduce a concept, which we call the $n$-filling property, for an action of a group on a separable unital $C^*$-algebra $A$.
If $A=C(\Om)$ is a commutative unital $C^*$-algebra and the action is induced by a group of homeomorphisms of $\Om$ then the $n$-filling property reduces to a weak version of hyperbolicity.
The $n$-filling property is used to prove that certain crossed product $C^*$-algebras are purely infinite and simple. A variety of group actions on boundaries of symmetric spaces and buildings have the $n$-filling property. An explicit example is the action of $\G=SL_n(\ZZ)$ on the projective $n$-space. 
\end{abstract}
\maketitle
\section*{Introduction}

Consider a $C^*$-dynamical system $(A,\alpha,\G)$ where $A$ is a separable unital $C^*$-algebra on which $\G$ acts by $*$-automorphisms.

\begin{definition}\label{ncnfilling}
Let $n \ge 2$ be a positive integer. We say that an action $\alpha : g \mapsto \alpha_g$ of $\G$ on $A$ is 
{\em $n$-filling} if, for all $b_1,b_2,\dots,b_n\in A^+$, with $\|b_j\|=1,\,1\le j \le n$, and for all $\e >0$, there exist $g_1,g_2,\dots,g_n\in\G$ such that $\sum_{j=1}^n \alpha_{g_j}(b_j)\ge 1-\e$. 
\end{definition}

If $A$ is a commutative unital $C^*$-algebra and $\alpha$ is induced by a group of homeomorphisms of the spectrum $\Om$ of $A$, then the $n$-filling property is equivalent to a generalized global version of hyperbolicity (Proposition \ref{comnfilling} below). In this setting, the definition was motivated by ideas from \cite {ad,ls} and \cite{bch}. The present article applies the $n$-filling property to give a proof that certain crossed product $C^*$-algebras are purely infinite and simple (Theorem \ref{ncnf}). In the commutative case, similar results were obtained in \cite{ad,ls} using local properties of the action. The paper \cite{ad} also considers more general groupoid $C^*$-algebras. Simple crossed product algebras have been constructed using the related concept of a strongly hyperbolic action in \cite[Appendix 2]{h}.

\begin{remark}
In order to prove the $n$-filling condition as stated in Definition \ref{ncnfilling}  it is sufficient to verify it for all $b_1,b_2,\dots,b_n$ in a dense subset $C$ of $A^+$. For then if $b_1,b_2,\dots,b_n\in A^+$, with $\|b_j\|=1,\,1\le j \le n$, and if $\e >0$, choose $c_1,c_2,\dots,c_n \in C$ such that $\|b_j-c_j\|<\frac{\e}{2n}$ for all $j$ and $\sum_{j=1}^n \alpha_{g_j}(c_j)\ge 1-\e/2$.  Write 
$$\sum_{j=1}^n \alpha_{g_j}(b_j-c_j)=x=x_+-x_-$$
where $x_+,x_-\in A^+$ and $x_+x_-=0$.
We have $x\ge -\e/2$ and therefore
$$\sum_{j=1}^n \alpha_{g_j}(b_j)=\sum_{j=1}^n \alpha_{g_j}(c_j)+x\ge1-\e/2-\e/2=1-\e.$$
\end{remark}

\bigskip
Suppose that $A=C(\Om)$, the algebra of continuous complex valued functions on a compact Hausdorff space $\Om$. If the action arises from an action of $\G$ on $\Om$ by homeomorphisms, then the $n$-filling condition can be expressed in the following way, which explains its name.

\begin{proposition}\label{comnfilling}
Let $\Om$ be an infinite compact Hausdorff space and let $\G$ be a group which acts on $\Om$ by homeomorphisms.
The induced action $\alpha$ of $\G$ on $C(\Om)$ is $n$-filling if and only if the following condition is satisfied:\ for any nonempty open subsets $U_1, \dots , U_n$ of $\Om$, there exist $g_1,\dots,g_n \in \G$ such that $g_1U_1 \cup \dots \cup g_nU_n = \Om$.
\end{proposition}

\begin{proof} If the action is $n$-filling, let $U_1, \dots , U_n$ be nonempty open subsets of $\Om$. There exist elements $b_1,b_2,\dots,b_n\in A^+$, with $\|b_j\|=1$, such that 
$\supp (b_j)\subset U_j$,\ $1\le j \le n$. By hypothesis there exist $g_1,g_2,\dots,g_n\in\G$ such that $\sum_{j=1}^n \alpha_{g_j}(b_j)\ge 1/2$.
Then if $\w\in\Om$ there exists $i\in\{1,2,\dots,n\}$ such that $\alpha_{g_i}(b_i)(\w)>0$. Therefore 
$g_i^{-1}\w\in U_i$,  i.e. $\w\in g_iU_i$. Thus $g_1U_1 \cup \dots \cup g_nU_n = \Om$.

Conversely, suppose the stated assertion holds. Fix $b_1,b_2,\dots,b_n\in A^+$, with $\|b_j\|=1,\,1\le j \le n$, and  let $\e >0$. For each $j$, the set  $U_j=\{\w\in\Om\,;\, b_j(\w)>1-\e\}$ is a nonempty and open.
Choose $g_1,\dots,g_n \in \G$ such that $g_1U_1 \cup \dots \cup g_nU_n = \Om$.
If $\w\in\Om$, then $g_i^{-1}\w\in U_i$ for some $i$ and so $\alpha_{g_i}(b_i)(\w)>1-\e$. Therefore $\sum_{j=1}^n \alpha_{g_j}(b_j)\ge 1-\e$. 
\end{proof}

\begin{remark}{\rm 
If the action of the group $\G$ on the space $\Om$ is topologically
transitive (in particular, if it is minimal) then the $n$-filling condition is equivalent to the following apparently weaker condition:
for each nonempty open subset $U$ of $\Om$, there exist $t_1,\dots,t_n \in \G$ such that
$t_1U \cup \dots \cup t_nU = \Om$.

In order to see this, suppose that  $U_1, \dots , U_n$ are nonempty open subsets
of $\Om$. There exists $g_2 \in \G$ such that $U_1 \cap g_2U_2 \ne \emptyset$.
Then there exists $g_3 \in \G$ such that $U_1 \cap g_2U_2 \cap g_3U_3 \ne \emptyset$.
Finally, there exists $g_n \in \G$ such that 
$U = U_1 \cap g_2U_2 \dots \cap g_nU_n \ne \emptyset$.
Then there exist $t_1,\dots,t_n \in \G$ such that
$t_1U \cup \dots \cup t_nU = \Om$ and so
$t_1U_1 \cup t_2g_2U_2 \dots \cup t_ng_nU_n = \Om$. 
}\end{remark}

\begin{definition}\label{n-filling}
Let $\phi(\G,\Om)$ be the smallest integer $n$ for which the conclusion of Proposition \ref{comnfilling} holds. Set $\phi(\G,\Om)= \infty$ if no such $n$ exists; that is, if the action is not $n$-filling for any integer $n$.
\end{definition}

Topologically conjugate actions have the same value of $\phi(\G,\Om)$.
It is easy to see that the notion of a $2$-filling action is equivalent to what is called a  {\it strong boundary action} in \cite{ls} and an {\it extremely proximal flow} in \cite{g}. The action of a word hyperbolic group on its Gromov boundary is $2$-filling \cite[Example 2.1]{ls}.
In our first example below (Example \ref{projectivespace}) we show that the canonical action of $\G=SL_n(\ZZ)$ on the projective space $\Pi=\PP^{n-1}(\RR)$ satisfies $\phi(\G,\Pi)=n$.

The final part of the paper is devoted to estimating $\phi(\G,\Om)$ for some group actions on the boundaries of affine buildings.
These estimates show that $\phi(\G,\Om)$ is not a stable isomorphism invariant for the algebra $C(\Om) \rtimes_r \G$ (Example \ref{notsii}).
\bigskip

\section{Purely infinite $C^*$-algebras from $n$-filling actions}\label{pisun}

\begin{definition}\label{properlyouter}
An automorphism $\alpha$ of a $C^*$-algebra $A$ is said to be {\em properly outer} if for each nonzero $\alpha$-invariant ideal $I$ of $A$ and for each inner automorphism $\beta$ of $I$ we have $\| \alpha |I - \beta \| =2$.
\end{definition}
We shall say that an action $\alpha : g \mapsto \alpha_g$ is  properly outer if for all $g\in \G\backslash\{e\}$, $\alpha_g$ is properly outer.

The purpose of this section is to prove the following result.

\begin{theorem}\label{ncnf}
Let $(A,\alpha,\G)$ be a $C^*$-dynamical system, where $A$ is a separable unital $C^*$-algebra. Suppose that for every nonzero projection $e\in A$ the hereditary $C^*$-subalgebra $eAe$ is infinite dimensional. Suppose also that the action $\alpha$ is $n$-filling and properly outer. Then the reduced crossed product algebra $B=A\rtimes_{\alpha,r}\G$ is a purely infinite simple $C^*$-algebra. 
\end{theorem}

\begin{remark}

\noindent  If $A=C(\Om)$, with $\Om$ a compact Hausdorff space, the condition that $eAe$ is infinite dimensional for every nonzero projection $e\in A$ says simply that the space $\Om$ has no isolated points.

It was shown in \cite[Proposition 1]{as} that if the action $\alpha$ is {\em topologically free} then  $\alpha$ is properly outer.
\end{remark}

\begin{proof} (Inspired by \cite[Theorem 5]{ls}.)
Denote by $E:B\to A$ the canonical conditional expectation. Fix $x\in B, \, x\ne 0$. In order to prove the result it is enough to show that there exist $y,z\in B$ such that $yxz=1$. Put $a=\frac{x^*x}{\|E(x^*x)\|}$. Let 
$0<\e<\frac{1}{2(2n+1)}$.
There exists $b\in C_c(\G,A)^+$ such that $\|a-b\|<\e$. Write $b=b_e+\sum_{g\in F}b_gu_g$, where $b_e=E(b)\ge 0$ and $F\subset \G\setminus\{e\}$ is finite.
Note that $\e>\|E(a-b)\|=\|E(a)-b_e\|\ge \big|1-\|b_e\|\big|$, and so $\| b_e\|^{-1}<1+2\e$.
It follows that 
$$\big\|a-\frac{b}{\| b_e\|}\big\|=\| b_e\|^{-1}\big\|(\| b_e\|-1)a +a-b\big\|
<(1+2\e)(\e\| a\|+\e)=\e(1+2\e)(1+\| a\|).$$
Choosing $b$ so that $\| a-b\|<\frac{\e}{3(1+\| a\|)}$ then replacing $b$ by 
$\frac{b}{\| b_e\|}$ shows that we can assume that $\| b_e\|=1$.

 Since $\alpha_g$ is properly outer for each $g\in F$, it follows from \cite[Lemma 7.1]{op} that there exists $y\in A^+$, $\| y\|=1$ such that
$\|b_e\|\ge\|yb_ey\|>\|b_e\|-\e/|F|$ and $\|yb_g\alpha_g(y)\|<\e/|F|$ for all $g\in F$.
Using Lemma \ref{ncnfm} below, we see that
there exists $c\in B$ such that $\|c\|\le\sqrt n$ and $c^*yb_eyc\ge1-3\e$.

Then 
\begin{equation*} \begin{split}
\| c^*yayc-c^*yb_eyc\|
&\le  \| c^*yayc-c^*ybyc\| + \| c^*ybyc-c^*yb_eyc\|  \\ 
&\le n\|a-b\| + n\|yby-yb_ey\| \\
&\le n\e +n\sum_{g\in F}
\| yb_gu_gyu_g^{-1}u_g\| \le 2n\e
\end{split}\end{equation*}

Therefore $c^*yayc$ is invertible since $(\|c^*yb_eyc)^{-1}\| \le \frac{1}{1-3\e}$ and
\begin{equation*}
\| 1-(c^*yb_eyc)^{-1}
(c^*yayc)\| 
\le \frac{2n\e}{1-3\e} < \frac{n}{2n-1} <1.
\end{equation*}
Setting $z=(c^*yayc)^{-1}$ we have 
$\|E(x^*x)\|^{-1}c^*yx^*\cdot x\cdot ycz=1.$
\end{proof}

It remains to prove Lemma \ref{ncnfm}. A preliminary observation is necessary.

\begin{lemma}\label{ncnfl}
Let $A$ be a unital $C^*$-algebra such that 
for every nonzero projection $e\in A$ the hereditary $C^*$-subalgebra $eAe$ is infinite dimensional. Let $b\in A^+,\,  \| b\|=1$ and let $\e>0$. For every integer $n\ge 1$ there exist elements  $b_1, b_2,\dots, b_n \in A^+$, with $\| b_j\|=1$, $bb_j=b_jb$, $\| bb_j\|\ge 1-\e$ and  $b_ib_j=0, for\ i\ne j$.
\end{lemma}

\begin{proof}
There are two cases to consider.

\noindent {\bf Case 1}. Suppose that $1$ is not an isolated point of $\Spec (b)$. Then there exist pairwise disjoint nonempty open sets  
$U_1, \dots , U_n$ contained in $\Spec (b) \cap [1-\e,1]$.
Let $C$ be the $C^*$-subalgebra of $A$ generated by $\{ b,1\}$.
By functional calculus, there exist $b_1, b_2,\dots, b_n \in C^+$, $\| b_j\|=1\,(1\le j \le n)$
with $\| bb_j\|\ge 1-\e$ and $b_ib_j=0,\, i\ne j$.

\noindent {\bf Case 2}. 
Suppose that $1$ is an isolated point of $\Spec (b)$. Then there exists a nonzero projection $e\in A$ such that
$be=eb= e$. By hypothesis the hereditary $C^*$-subalgebra $eAe$ is infinite dimensional. Therefore every masa of $eAe$ is infinite dimensional \cite[p. 288]{kr}. Inside such an infinite dimensional masa of $eAe$ we can find positive elements $b_1, b_2,\dots, b_n$, $\| b_j\|=1\,(1\le j \le n)$
with  $b_ib_j=0,\, i\ne j$. Then  $bb_j=b(eb_j)=eb_j=b_j=b_jb$ and $\|bb_j\|=\|b_j\|=1$
for  $1\le j\le n$.
\end{proof}

\begin{lemma}\label{ncnfm}  Let $(A,\alpha,\G)$ be as in the statement of Theorem \ref{ncnf}, let $0<\e<1/3$ and let $b\in A^+$, with $1-\e\le \| b\|\le1$. Then there exists $c\in B$ such that $\| c \|\le \sqrt n$ and $c^*bc\ge 1-3\e$.
\end{lemma}

\begin{proof}
By Lemma \ref{ncnfl}, there exist $b_1, b_2,\dots, b_n \in A^+$, with $\| b_j\|=1$, $bb_j=b_jb$, $b_ib_j=0$ for $i\ne j$, and  $\| bb_j\|\ge 1-2\e$.
Since the action is $n$-filling, there exist $g_1,g_2,\dots,g_n\in\G$ such that $\sum_{i=1}^n 
\frac{1}{\| bb_i\| }\alpha_{g_i}(bb_i)\ge 1-\e$. Therefore $\sum_{i=1}^n\alpha_{g_i}(bb_i)\ge (1-\e)(1-2\e)\ge 1-3\e.$ 
Put $c=\sum_{j=1}^n{\sqrt b_j}u_{g_j}^{-1} \in B$.

Now $c^*c=\sum_{i,j}u_{g_i}{\sqrt b_i}{\sqrt b_j}u_{g_j}^{-1}=\sum_{i=1}^n \alpha_{g_i}(b_i)\le n$ and so $\| c\|\le\sqrt n$.
Finally, we have
$c^*bc= \sum_{i,j}u_{g_i}{\sqrt b_i}b{\sqrt b_j}u_{g_j}^{-1}= \sum_{i=1}^n \alpha_{g_i}(bb_i)   \ge 1-3\e.$ 
 \end{proof}

\bigskip
\section{examples}\label{examples}
We now give some explicit examples of $n$-filling actions.

\begin{example}\label{projectivespace}
For the canonical action of $\G=SL_n(\ZZ)$ on the projective space $\Pi=\PP^{n-1}(\RR)$, we have $\phi(\G,\Pi)=n$.
\end{example}

\begin{proof} 
Denote by $u \mapsto [u]$ the canonical map from $\RR^n$ onto $\Pi$.

We first show that the action of $\G$ on $\Pi$ is not $(n-1)$-filling. Choose a linear subspace $E$ of $\RR^n$ of dimension $n-1$. Let $U=\Pi\setminus [E]$,
which is a nonempty open subset of $\Pi$. If $t_j \in \G$ ($1\le j \le n-1$) then $t_1U \cup \dots \cup t_{n-1}U \ne \Pi$. For the subspace $t_1E \cap \dots \cap t_{n-1}E$ of $\RR^n$ has dimension at least one, and so contains a nonzero vector $v$. Then $[v] \notin \bigcup_{j=1}^{n-1}t_jU$. Thus the action$(\G,\Pi)$ is not $(n-1)$-filling. It remains to show that it is $n$-filling. For this we use ideas from \cite[Example 1]{bch}.

We claim that there exists a basis $\{u_1,u_2,\dots,u_n\}$ for $\RR^n$, elements $g_1,g_2,\dots,g_n \in \G$, and (compact) sets $K_1,K_2,\dots,K_n \subset \Pi$ with $K_1\cup K_2\cup\dots\cup K_n=\Pi$, and with the following property: for any open neighbourhood $U_j$ of $[u_j]$ ($1\le j \le n$) there exists a positive integer $N_j$ such that $g_j^nK_j \subset U_j$ for all $n\ge N_j$. It follows that the action is $n$-filling.  For let $U_1, \dots , U_n$ be nonempty open subsets
of $\Pi$. Since the action of $\G$ on $\Pi$ is minimal, we may assume that $[u_j]\in U_j$ ($1\le j \le n$). Let $t_j=g_j^{-N_j}$, so that $K_j \subset t_jU_j$ ($1\le j \le n$). Then $t_1U_1 \cup \dots \cup t_nU_n = \Pi$.

It remains to verify our claim. Fix a positive integer $k \ge 4$ and let $a=\frac{2}{\sqrt{k^2+4k}+k}$, $b=\frac{\sqrt{k^2+4k}-k}{2}$. Consider the matrices
$A = \bigl( \begin{smallmatrix}
k+1&k\\ 1&1
\end{smallmatrix} \bigr)$
and 
$B = \bigl( \begin{smallmatrix}
1&1\\ k&k+1
\end{smallmatrix} \bigr)$
in $SL_2(\ZZ)$. These matrices have eigenvalues $\lambda_+ = 1+\frac{1}{a}$, $\lambda_- = 1-b$, which satisfy $0<\lambda_-<1<\lambda+$. The corresponding eigenvectors for $A$ are 
$\bigl( \begin{smallmatrix} 1\\ a \end{smallmatrix} \bigr)$
and $\bigl( \begin{smallmatrix} -b\\ 1\end{smallmatrix} \bigr)$;
for $B$ they are 
$\bigl( \begin{smallmatrix} a\\ 1 \end{smallmatrix} \bigr)$
and $\bigl( \begin{smallmatrix} 1\\ -b\end{smallmatrix} \bigr)$.
If $1\le j\le n-1$ let 
\begin{equation*}
\begin{aligned}
g_j=
\begin{pmatrix}1&0&\dots&0\\
               0&1&\dots&0\\
                & &     & \\
                & & A   & \\
                & &     & \\
               0&0&\dots&1
\end{pmatrix}
\end{aligned}
\quad\text{,}\quad
\begin{aligned}
u_j=
\begin{pmatrix}0\\
               0\\
                \\
                1\\
                a\\
                 \\
                0
\end{pmatrix}
\end{aligned}
\quad\text{,}\quad
\begin{aligned}
v_j=
\begin{pmatrix}0\\
               0\\
                \\
                -b\\
                1\\
                 \\
                0
\end{pmatrix}
\end{aligned}
\end{equation*}
where A occupies the $j$ and $j+1$ rows and columns and the nonzero entries of the vectors are in rows $j$ and $j+1$. Also let 
\begin{equation*}
\begin{aligned}
g_n=
\begin{pmatrix}1&0&\dots&0\\
               0&1&\dots&0\\
                & &     & \\
                & &   & \\
               &&     & \\
               &&     &B
\end{pmatrix}
\end{aligned}
\quad\text{,}\quad
\begin{aligned}
u_n=
\begin{pmatrix}0\\
               0\\
                \\
                \\
                \\
                a\\
                1
\end{pmatrix}
\end{aligned}
\quad\text{,}\quad
\begin{aligned}
v_n=
\begin{pmatrix}0\\
               0\\
                \\
                \\
                \\
                1 \\
                -b
\end{pmatrix} .
\end{aligned}
\end{equation*}

Let $R={\rm max}(\frac{1+a}{1-b},\frac{1+ab}{1-b})=\frac{1+a}{1-b}$.
For $1\le j\le n-1$ let
$$K_j=\{[\xi_ju_j+\eta_jv_j+\sum_{l\ne j,j+1}\xi_le_l]\ ;\ \xi_j\ne 0,
\big|\frac{\eta_j}{\xi_j}\big|\le R, \big|\frac{\xi_l}{\xi_j}\big|\le R, l\ne j,j+1\},$$
$$K_n=\{[\xi_nu_n+\eta_nv_n+\sum_{l\ne n-1,n}\xi_le_l]\ ;\ \xi_n\ne 0,
\big|\frac{\eta_n}{\xi_n}\big|\le R, \big|\frac{\xi_l}{\xi_n}\big|\le R, l\ne n-1,n\}.$$

Direct computation shows if $[x] \in \Pi$ then $[x]\in K_j$, where $|x_j| = {\rm max}_{1\le l \le n}|x_l|$. Therefore $\Pi = \bigcup_{j=1}^nK_j$.

Let $\epsilon >0$ and consider the basic open neighborhood $U_j$ of $[u_j]$ defined by
$$U_j=\{[\xi_ju_j+\eta_jv_j+\sum_{l\ne j,j+1}\xi_le_l]; \xi_j\ne 0,
\big|\frac{\eta_j}{\xi_j}\big| < \epsilon, \big|\frac{\xi_l}{\xi_j}\big| < \epsilon, l\ne j,j+1\}.$$
Let $N > \frac{{\rm log}(R/\epsilon)}{{\rm log}(\lambda_+)}$. Recall that $0<\lambda_-<1<\lambda_+$. Therefore $\frac{R}{\lambda_+^N}<\epsilon$. 

For $m \ge N$ and $[\xi_ju_j+\eta_jv_j+\sum_{l\ne j,j+1}\xi_le_l]\in K_j$, we have
$$
g^m[\xi_ju_j+\eta_jv_j+\sum_{l\ne j,j+1}\xi_le_l]
=[\lambda_+^m\xi_ju_j+\lambda_-^m\eta_jv_j+\sum_{l\ne j,j+1}\xi_le_l].$$
Now $\big|\frac{\lambda_-^m\eta_j}{\lambda_+^m\xi_j}\big|\le 
\frac{1}{\lambda_+^m}\big|\frac{\eta_j}{\xi_j}\big|\le
\frac{R}{\lambda_+^m}< \epsilon$, and for $l\ne j,j+1$,
$\big|\frac{\xi_l}{\lambda_+^m\xi_j}\big|\le 
\frac{1}{\lambda_+^m}\big|\frac{\xi_l}{\xi_j}\big|\le
\frac{R}{\lambda_+^m}< \epsilon$.   

This means that $g_j^mK_j \subset U_j$ for all $m\ge N$.
\end{proof}

\begin{remark}{\rm
The fact that the action of $SL_3(\ZZ)$ on the projective plane $\PP^{2}(\RR)$ is not $2$-filling can also be seen in a different way. More generally the action of a group $\G$ on a non-orientable compact surface $\Om$ cannot be $2$-filling.
For let $M$ be a closed subset of $\Om$ homeomorphic to a M\"obius band, let $U_1=M^c$ and let $U_2 \subset \Om$ be homeomorphic to an open disc in $\RR^2$. Then it is impossible to have $t_1U_1 \cup t_2U_2 = \Om$ for $t_1,t_2 \in \G$. For  $t_2^{-1}t_1(M)$ would be a homeomorphic copy of a M\"obius band embedded in the disc $U_2$. To see that this is impossible note that a Mobius band is not disconnected by its centre circle, and apply the Jordan curve theorem. }\end{remark}

\begin{definition}\label{attracting}
Let the group $\G$ act on the topological space $\Om$. An element $g \in \G$ is said to have an attracting fixed point $x\in \Om$ if $gx=x$ and there exists a neighbourhood $V_x$ of $x$ such that ${\displaystyle\lim_{n \to \infty}}g^n(V_x)=\{x\}$.
\end{definition}

\begin{remark}{\rm
Let $G$ be a noncompact semisimple real algebraic group and let $\G$ be a Zariski-dense subgroup of $G$. Consider the action of $G$ on its Furstenberg boundary $G/P$, where $P$ is a minimal parabolic subgroup of $G$. It follows from \cite[Appendice]{bel} that there exist elements $g\in \G$ which have attracting fixed points in $G/P$. In fact the set $H$ of all such elements $g \in \G$ is Zariski-dense in $G$: the elements of $H$ are called {\it h-regular} in \cite{bel} and {\it maximally hyperbolic} in \cite{bch}.

It follows from a result of {\sc H. Furstenberg} \cite[Theorem 5.5, Corollary]{fur} that if $G$ is a semisimple group with finite centre which acts minimally on a locally compact Hausdorff space $\Om$ with an attracting fixed point, then $\Om$ is necessarily a compact homogeneous space of $G$. 
}
\end{remark}

The following result shows that many of the actions considered in \cite{ad,ls} are $n$-filling for some integer $n$.

\begin{proposition}\label{minattract}
Let $\Om$ be a compact Hausdorff space and let $(\Om,\G)$ be a minimal action. Suppose that there exists an element $g\in \G$ which has an attracting fixed point in
$\Om$. Then the action $(\Om,\G)$ is $n$-filling for some integer $n$.
\end{proposition}

\begin{proof}
Choose $x\in \Om$ with $gx=x$ and an open neighbourhood $V_x$ of $x$ such that $ \lim_{n \to \infty}g^n(V_x)=\{x\}$. Since the action is minimal, the family $\{hV_x ; h \in \G\}$ forms an open covering of $\Om$. By compactness, there exists a finite subcovering $\{h_1V_x,h_2V_x,\dots,
h_nV_x\}$. 

Let $U_1, \dots , U_n$ be nonempty open subsets of $\Om$. Since the action of $\G$ on $\Om$ is minimal, we may choose elements $s_j \in \G$ such that  $h_jx\in s_jU_j$ ($1\le j \le n$). For $1\le j \le n$, choose an integer $N_j$ such that $g^{N_j}V_x \subset h_j^{-1}s_jU_j$.  Then $h_jV_x \subset t_jU_j$, where $t_j=h_jg^{-N_j}h_j^{-1}s_j$. Therefore $t_1U_1 \cup \dots \cup t_nU_n = \Om$.
\end{proof}

\begin{remark}\label{bech}{\rm
Consider the action of a noncompact semisimple real algebraic group $G$ on its Furstenberg boundary $G/P$. Let $\G$ be a Zariski-dense subgroup of $G$ and let $n(W)$ be the order of the Weyl group. 
In this case one can be more precise: the action $(G/P,\G)$ is $n(W)$-filling. The proof follows from the remarks in \cite[page 127]{bch}.  In the next section we prove an analogue of this result for groups acting on affine buildings.
}
\end{remark}

Recall that an action $(\Om^{\prime},\G)$ is said to be a {\it factor} of the action $(\Om,\G)$ if there is a continuous equivariant surjection from $\Om$ onto $\Om^{\prime}$.

\begin{proposition}\label{factor}
Suppose that the action $(\Om,\G)$ is $n$-filling and that $(\Om^{\prime},\G)$ is a factor of $(\Om,\G)$. Then $(\Om^{\prime},\G)$ is an $n$-filling action.
\end{proposition}

\begin{proof}
This is an easy consequence of the definitions.
\end{proof}
\bigskip
\section{Group actions on boundaries of affine buildings}

We now turn to some examples which motivated our definition of an $n$-filling action. They are discrete analogues of those referred to Remark \ref{bech}.
We show that if a group $\G$ acts properly and cocompactly on
an affine building $\D$ with boundary $\Om$, then the induced action on $\Om$ is a $n$-filling, where $n$ is the number of boundary points of an apartment in $\D$. If $\D$ is the affine Bruhat-Tits building of a linear group then $n$ is the order of the associated spherical Weyl group.

  An {\sl apartment} in $\D$ is a subcomplex of $\triangle$ isomorphic to an affine Coxeter complex. Each apartment inherits a natural metric from the Coxeter complex, which gives rise to a well-defined metric on the whole building \cite[Chapter IV.3]{bro}. Every geodesic of $\D$ is a straight line in some apartment. A {\sl sector} (or {\sl Weyl chamber}) is a sector based at a special vertex in some apartment \cite{ron}. Two sectors are  {\sl equivalent} (or parallel) if their intersection contains a sector. The boundary $\Omega$ is defined to be the set of equivalence classes of sectors in $\triangle$.  Fix a special vertex $x$.  For any $\omega \in \Omega$ there is a unique sector $[x,\omega)$ in the class $\omega$ having base vertex $x$ \cite[Theorem 9.6, Lemma 9.7]{ron}. In the terminology of \cite[Chapter VI.9]{bro}
$\Omega$ is the set of chambers of the building at infinity $\triangle^\infty$. 
Topologically, $\Omega$ is a totally disconnected compact Hausdorff space and a basis for the topology is given by sets of the form
$$
\Omega_x(v) = \left \{ \omega \in \Omega : [x,\omega) \hbox { contains } v \right \}
$$
where $v$ is a vertex of $\triangle$. See \cite[\S 2]{cms} for the  $\widetilde A_2$ case, which generalizes directly.

 We will need to use the fact that $\Om$ also has the structure of a spherical building \cite[Theorem 9.6]{ron}, and its apartments are topological spheres. 

\begin{definition}\label{opposite}
 Two boundary points $\w, \wa$ in $\Om$ are said to be {\it opposite} \cite[IV.5]{bro} if the distance between them is the diameter of the spherical building $\Om$. Opposite boundary points are opposite in a spherical apartment of $\Om$ which contains them; this apartment is necessarily unique.  Two subsets of $\Om$ are opposite if each point in one set is opposite each point in the other.   

We define $\cO(\w)$ to be the set of all $\w' \in \Om$ such that  $\w'$ is opposite to $\w$. It is easy to see that $\cO(\w)$ is an open set.
\end{definition}

\begin{lemma}\label{op apartment} If $\w \in \Om$ and $\cA$ is an 
apartment in $\D$, then there exists a boundary point $\wa$ of $\cA$
such that $\wa$ is opposite $\w$.
\end{lemma}

{\sc Proof:} 
Consider the geometric realization of the spherical building $\Om$. By \cite[Theorem (A.19)]{ron}, the subcomplex $\Om^{\prime}$ obtained from $\Om$ by deleting all chambers opposite $\w$ is geodesically contractible. However this is impossible if $\Om^{\prime}$ contains the spherical apartment of $\Om$ made up of the boundary points of $\cA$. 
\qed

\begin{corollary}\label{six} If $\w_1, \dots, \w_n$ are the  boundary points of an apartment then $$\Om= \cO(\w_1) \cup \dots \cup \cO(\w_n)$$
\end{corollary}

\begin{remark}
{\rm  The union is not disjoint in general, as is seen by considering the example of a tree.}
\end{remark}

\begin{lemma}\label{op} Two chambers $\w_1, \w_2$ in $\Om$ are opposite if and only if they are represented by opposite sectors $S_1, S_2$ with the same base vertex in some apartment of $\triangle$. Moreover if two sectors $S_1, S_2$ in an apartment $\cA$ with the same base vertex represent opposite elements $\w_1, \w_2$ in $\Om$, then $S_1, S_2$ are opposite sectors and $\cA$ is the unique apartment containing them. 
\end{lemma}

{\sc Proof:}
Suppose that $\w_1, \w_2$ in $\Om$ are opposite. There exists an apartment $\cA$ containing sectors $S_1, S_2$ representing  $\w_1, \w_2$ respectively \cite[Proposition 9.5]{ron} or \cite[VI.8,Theorem]{bro}. By taking parallel sectors, we may assume that $S_1, S_2$  have the same base vertex $x \in \cA$.
The sectors of $\cA$ based at $x$ correspond to the chambers of an apartment in $\Om$ containing $\w_1, \w_2$ \cite[Theorem 9.8]{ron}. Therefore $S_1, S_2$ are opposite
sectors. The converse is clear.

The final assertion follows from \cite[VI.9, Lemma 2 and IV.5 Theorem 1]{bro}.
\qed

\begin{remark}\label{opposite triangles}{\rm
\item (a) It is not necessarily true that if $\w_1, \w_2$ in $\Om$ are opposite then the sectors $[z,\w_1)$, $[z,\w_2)$ based at any vertex $z$ are opposite sectors in some apartment. 
\item (b) If $C_1, C_2$ are opposite chambers with a common  vertex $x$ in an apartment, then $\Om_x(C_1)$ and $\Om_x(C_2)$ are opposite sets in $\Om$.
}\end{remark} 
\bigskip

Suppose that a group $\G$ acts properly and cocompactly on an affine building $\D$ of dimension $n$. An apartment $\cA$ in $\D$ is said to be {\it periodic} if there is a subgroup $\G_0 < \G$ preserving $\cA$ such that $\G_0 \backslash \cA$
is compact \cite[6.$B_3$]{gro}. Note that $\G_0$ is commensurable with $\ZZ^n$, and this concept coincides with the notion of periodicity described in \cite{mz},\cite{rr} for buildings of type $\tilde A_2$.  In \cite{bb}, a periodic apartment is called $\G$-closed. This terminology makes it clear that periodicity
depends upon the choice of the group $\G$ acting on the building.

It is important to observe that there are many periodic apartments. In fact, according to \cite[Theorem 8.9]{bb}, any compact subset of an apartment is contained in some periodic apartment 

Now let $\cA_0$ be a periodic apartment, and fix a special vertex $z$ in $\cA_0$.
Choose a pair of opposite sectors $W^+$, $W^-$ in $\cA_0$ based at $z$.
Denote by $\w^{\pm}$ the boundary points represented by $W^{\pm}$, respectively.
By periodicity of the apartment there is a periodic direction represented by a line $L$ in any of the sector directions of $\cA_0$. For definiteness choose this direction to be that of the sector $W^+$. This means that there is an element $u \in \G$ which leaves $L$ invariant and translates the apartment $\cA_0$ in the direction of $L$. (In the terminology of \cite{bb, moz}, $L$ is said to be an {\it axis} of $u$.) Then $u^n\w^+ = \w^+$, $u^n\w^- = \w^-$ for all $n \in \ZZ$.
Moreover  $u^nz$ is in the interior of $W^+$ for $n > 0$ and in the interior of $W^-$ for $n < 0$. (See Figure \ref{W+}. Here and in what follows, the figures illustrate the case of a building $\D$ of type $\tilde A_2$, where each apartment contains precisely six sectors based at a given vertex.) The element $u$ above is the analogue of the maximally hyperbolic elements in \cite{bch}.

\refstepcounter{picture}
\begin{figure}[htbp]\label{W+}
{}\hfill
\beginpicture
\setcoordinatesystem units <0.5cm,0.866cm>    
\setplotarea x from -6 to 6, y from -4 to 4         
\put {$\bullet$} at 0 0
\put {$\bullet$} at 0 3
\put {$\bullet$} at 0 -3
\put {$z$}             [b l]  at    0.3   0.1
\put {$u^nz$}  [l]  at  0.2  3  
\put {$u^{-n}z$}  [l]  at  0.2  -3    
\put {$W^+$}   at    0  2
\put {$W^-$}   at    0  -2
\putrule from -6    0     to  6    0
\setlinear
\plot -3  3   3 -3 /
\plot -3 -3   3  3 /
\setdashes
\putrule from 0 -1.5  to  0 1.5
\putrule from 0 2.5  to  0 3.5
\putrule from 0 -2.5  to  0 -3.5
\endpicture
\hfill{}
\caption{The periodic apartment $\cA_0$.}
\end{figure}

The following crucial result shows that $\w^-$ is an attracting fixed point for $u^{-1}$.

\begin{proposition}\label{claim} 
Let $\cA_0$ be a periodic apartment and choose a pair of opposite boundary points $\w^{\pm}$. Let $u \in \G$ be an element which translates the apartment $\cA_0$ in the direction of
$\w^+$. Then $u^{-1}$ attracts $\cO(\w^+)$ towards $\w^-$, that is: for each compact subset $G$ of $\cO(\w^+)$ we have  $\lim_{n \to \infty} u^{-n}(G) = \{\w^-\}$.
\end{proposition}

{\sc Proof:} We use the notation introduced above. Let $\w \in \cO(\w^+)$.
By considering a retraction of $\triangle$ centered at $\w^+$
\cite[p.170, VI.8, Theorem]{bro}, we see that $\triangle$ is a union of apartments
which contain a subsector of $W^+$. Moreover for any sector $W$ representing
$\w$ there are subsectors $V^+ \subset W^+$ and $V \subset W$ which lie
in a common apartment $\cA$. Replacing $V^+$ by a subsector, we may
assume that $V^+$ has base vertex $u^Nz$ for some $N$, that is $V^+=[u^Nz,\w^+)$.
Replacing $V$ by a parallel sector in $\cA$ we may also assume that $V$ has base vertex $u^Nz$. By Lemma \ref{op}, $V$ lies in the apartment $\cA$ as shown in Figure \ref{apartment A}.

\refstepcounter{picture}
\begin{figure}[htbp]\label{apartment A}
{}\hfill
\beginpicture
\setcoordinatesystem units <0.5cm,0.866cm>    
\setplotarea x from -6 to 6, y from -4 to 4         
\put {$\bullet$} at 0 0
\put {$u^Nz$}             [b l]  at    0.35   0.1
\put {$V^+$}   at    0  2
\put {$[{u^Nz},\w^+)$}   at    0  3.5
\put {$[{u^Nz},\w)$}  at     0 -3.5
\put {$V$}  at    0 -2
\putrule from -6    0     to  6    0
\setlinear
\plot -3  3   3 -3 /
\plot -3 -3   3  3 /
\endpicture
\hfill{}
\caption{The apartment $\cA$.}
\end{figure}

For each $N \ge 0$ let $G_N$ denote the set of all boundary points $\w \in \cO(\w^+)$ such that $[{u^Nz},\w)$ and $[{u^Nz},\w^+)$ are opposite sectors in some apartment $\cA^{(N)}$. Then $G_0 \subset G_1 \subset G_2 \subset \dots$ is an increasing family of compact open sets and we have observed above that $\bigcup_{N=0}^{\infty}G_N =\cO(\w^+)$.
The result will follow if we show that $\lim_{n \to \infty} u^{-n}(G_N) = \{\w^-\}$ for each $N \ge 0$. It is clearly enough to consider the case $N=0$.

\begin{figure}[htbp]
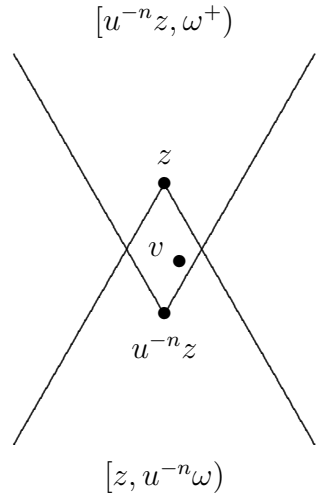
\label{rabbit}
{}\hfill
\beginpicture
\setcoordinatesystem units <0.5cm,0.866cm>    
\setplotarea x from 0 to 9, y from -4 to 4         
\put {$\bullet$} at  5.4 -0.2
\put {$v$}  [r]  at 5.0 0
\put {$\bullet$} at 5 -1
\put {$u^{-n}z$} [t]  at 5 -1.3
\put {$\bullet$} at 5  1
\put {$z$} [b]  at 5  1.3
\put {$[{u^{-n}z},\w^+)$}   at    5  3.5
\put {$[{z},u^{-n}\w)$}  at     5 -3.5
\setlinear
\plot 5 -1  9 3    /
\plot 5 -1   1 3 /
\plot 5 1  9 -3    /
\plot 5 1   1 -3 /
\endpicture
\hfill{}
\caption{Sectors in the apartment $u^{-n}\cA^{(0)}$.}
\end{figure}

Consider a basic open neighbourhood of $\w^-$ of the form $\Om_z(v)$, where $v \in [z,\w^-)\subset \cA_0$. Choose an integer $p \ge 0$ such that $u^nv \in [z,\w^+)$ for all $n \ge p$. If $\w \in G_0$ then $u^nv \in [u^nz,\w)$ (that is $v \in [z,u^{-n}\w)$) for all $n \ge p$. (See Figure \ref{rabbit}.)
This means that $u^{-n}\w \in \Om_z(v)$ for all $n \ge p$. Thus $u^{-n}(G_0)\subset \Om_z(v)$ for all $n \ge p$. This proves the result.
\qed        

\refstepcounter{picture}

\begin{theorem}\label{6strong}Suppose that a group $\G$ acts properly and cocompactly on the vertices of  an affine building $\D$ with boundary $\Om$. Let $k$ denote the number of boundary points of an apartment of $\D$. Then the action $(\Om,\G)$ is $k$-filling.
\end{theorem}

{\sc Proof:} Let $U_1, \dots , U_k$ be nonempty open subsets of $\Om$. Let $\cA_0$ be a periodic apartment with boundary points $\w_j, 1 \le j \le k$. By minimality of the action we can assume that $\w_j \in U_j, 1 \le j \le k$. By Corollary \ref{six}, we have $\Om = \cO(\w_1) \cup \dots \cup \cO(\w_n)$. It follows from the existence of a partition of unity that there exist compact sets $K_j \subset \cO(\w_j)$, $1 \le j \le k$ such that $\Om = K_1 \cup \dots \cup K_k$.

Let $u_j \in \G$ translate the apartment $\cA_0$ in the direction of $\w_j$, $1 \le j \le k$. Then by Proposition \ref{claim}, there exists $N_j \ge 0$ such that $u_j^{-n}K_j \subset U_j$ whenever $n \ge N_j$, $1 \le j \le k$.
In other words, $K_j \subset u_j^nU_j$ whenever $n \ge N_j$, $1 \le j \le k$. Let
$t_j=u_j^{N_j}$. Then  $$\Om = K_1 \cup \dots \cup K_k \subset t_1U_1 \cup \dots \cup t_kU_k$$ as required.
\qed 

\begin{remark}\label{not2filling}
The action of an $\tA_2$ group $\G$ on the boundary $\Om$ of the associated building is $6$-filling. We do not know the precise value of $\phi(\G,\Om)$, but it is certainly greater than $2$.  
To see this, fix a point $\w_0 \in \Om$ and choose $U$
to be a nonempty open set opposite $\w_0$. If $t_1,t_2 \in \G$
then $t_1U$ and $t_2U$ are opposite the boundary points
$t_1\w_0$ and $t_2\w_0$ respectively and therefore cannot
cover $\Om$. To see this, choose a hexagonal apartment
of $\Om$ which contains $t_1\w_0$ and $t_2\w_0$ and choose
a chamber $\wa$ in this apartment which is not opposite 
$t_1\w_0$ or $t_2\w_0$. Then $\wa$ cannot lie in
$t_1U \cup t_2U$.
Therefore $2< \phi(\G,\Om) \le 6$.
\end{remark}

\bigskip
\section{Purely infinite simple $C^*$-algebras }

Throughout this section we consider only affine buildings of type $\tA_2$.
The $\widetilde A_2$ buildings are a particularly natural setting for our investigation. They are the simplest two-dimensional buildings, but they do not necessarily arise from linear groups. Crossed product $C^*$-algebras associated with them have been studied in \cite{rs,rs'}.  
In this case the building $\triangle$ is a simplicial complex whose maximal simplices ({\sl chambers}) are triangles.  An apartment of $\triangle$ is a subcomplex isomorphic to the Euclidean plane tessellated by equilateral triangles.  
 
The boundary $\Om$ may be identified with the flag complex of a projective plane $(P,L)$ \cite[page 81]{bro}. Flags will be denoted $(x_1,x_2)$ where $x_1 \in x_2$. If we identify chambers of $\Om$ with sectors based at a fixed vertex $v_0$ of type $0$, then a sector wall whose base panel is of type $1$ corresponds to an element of $P$ and  a sector wall whose base panel is of type $2$ corresponds to an element of $L$ \cite[Section 9.3]{ron}.
$P$ is the {\it minimal boundary} of $\triangle$ and has been studied in \cite{cms}, where it is denoted $\Om^l$. The topology on $P$ comes from the natural quotient map $\Om \to P$. Moreover the action of $\G$ on $\Om$ induces an action on $P$. Similar statements apply to $L$, and there is a homeomorphism $P \cong L$. 

\refstepcounter{picture}
\begin{figure}[htbp]
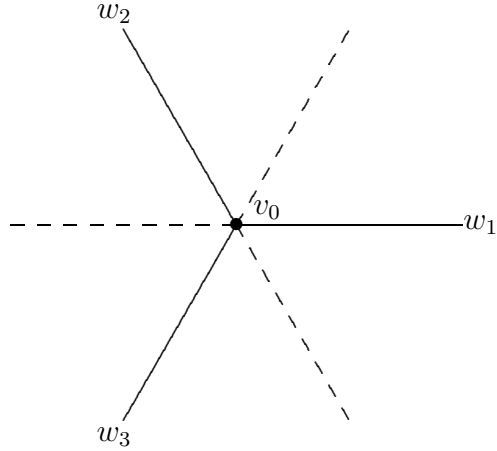
\label{P/L}
{}\hfill
\beginpicture
\setcoordinatesystem units <0.5cm,0.866cm>    
\setplotarea x from -6 to 6, y from -4 to 4         
\put {$\bullet$} at 0 0
\put {$v_0$}             [b l]  at    0.45   0.1
\put {$w_1$}               at    6.5 0
\put {$w_2$}               at    -3.25 3.25
\put {$w_3$}               at    -3.25 -3.25
\putrule from 0    0     to  6    0
\setlinear
\plot -3  3   0 0 /
\plot -3 -3   0 0 /
\setdashes
\putrule from -6    0     to  0   0
\plot 0 0   3 -3 /
\plot 0 0   3  3 /
\endpicture
\hfill{}
\caption{Sector walls $w_1$,$w_2$,$w_3$ corresponding to points in $P$.}
\end{figure}

From now on assume that the group $\G$ is an $\widetilde A_2$ group: that is $\G$ acts simply transitively in a type rotating manner on the vertices of an affine building $\D$ of type $\widetilde A_2$.

\begin{proposition}\label{topfree} The actions $(\Om,\G)$, $(P,\G)$ are topologically free. That is, if $g \in \G \backslash \{e\}$ then
$${\rm Int}\{ \w \in \Om : g\w = \w \} = \emptyset$$
$${\rm Int}\{ w \in P : gw = w \} = \emptyset$$
\end{proposition}

{\sc Proof:} The statement for the action on $\Om$ is proved in \cite[Theorem 4.3.2]{rs}.

\begin{figure}[htbp]
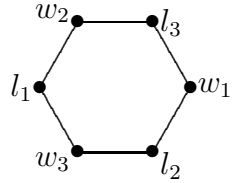
\label{P-L apartment}
{}\hfill
\beginpicture
\setcoordinatesystem units <0.5cm,0.866cm>    
\setplotarea  x from -2.5 to 2.5,  y from -1.5 to 1.5
\put {$\bullet$} at 1 1
\put {$\bullet$} at 1 -1
\put {$\bullet$} at -1 1
\put {$\bullet$} at -1 -1
\put {$\bullet$} at 2 0
\put {$\bullet$} at -2 0
\put {$l_3$} [l] at  1.2 1
\put {$l_2$} [t,l] at  1.2 -1
\put {$w_3$} [r,t] at  -1.2 -1
\put {$w_2$} [b,r] at  -1.2  1
\put {$w_1$} [l] at 2.2 0
\put {$l_1$} [r] at -2.2 0

\putrule from -1 1 to 1 1
\putrule from -1 -1 to 1 -1
\setlinear \plot 1 1  2 0  1 -1   /
\setlinear \plot -1 1  -2 0  -1 -1   /
\endpicture
\hfill{}
\caption{The apartment $\mathcal A_0$}
\end{figure}

Suppose that the result fails for the action on $P$.
Then there exists an open set $V \subset P$ such that
$gw = w$ for all $w \in P$.
Let $\widetilde V = \pi^{-1}(V)$, where
$\pi: \Om \to P$ is the quotient map.
Then $\widetilde V$ is a nonempty open subset of $\Om$.
By \cite[Proposition 4.3.1]{rs}, $\widetilde V$ contains
all six boundary points of some apartment $\mathcal A$ of $\triangle$.
These boundary points are the six chambers of an apartment
$\mathcal A_0$ in $\Om$, as illustrated in Figure \ref{P-L apartment}.
 The apartment $\mathcal A_0$ contains
three points $w_1, w_2, w_3 \in P$. These three points
lie in $V$ and hence are fixed by $g$. It follows that
the lines $l_1, l_2, l_3 \in L$ are also fixed by $g$.
Therefore each boundary point of $\mathcal A_0$ is fixed by $g$.
By the proof of \cite[Theorem 4.3.2]{rs}, it follows that 
$g\mathcal A = \mathcal A$ and $g$ acts by translation on $\mathcal A$.
The same is true for all nearby apartments $\mathcal A'$, since
the corresponding walls $w_1', w_2', w_3' \in P$
will also be fixed by $g$, if they belong to $V$.
The argument of \cite[Theorem 4.3.2]{rs} now
gives a contradiction.
\qed 

\refstepcounter{picture}

\begin{proposition}\label{main} 
If $\G$ is an $\widetilde A_2$ group, then the algebras $C(\Om) \rtimes \G$, $C(P) \rtimes \G$  are simple purely infinite $C^*$-algebras.
\end{proposition}

{\sc Proof:} The actions are topologically free by Proposition \ref{topfree} and hence properly outer \cite[Proposition 1]{as}. Moreover they are 
$6$-filling by Theorem \ref{6strong}.
The result follows from Theorem \ref{ncnf}.
\qed

We now give examples of properly outer actions $(\Om_i,\G_i)$, $i=1,2$, with $\phi(\G_1,\Om_1)=2$ and 
$\phi(\G_2,\Om_2)>2$ but for which $C(\Om_1) \rtimes \G_1$ is stably isomorphic to 
$C(\Om_2) \rtimes \G_2$.

\begin{example}\label{notsii}
{\rm
Let $\G_1 \subset \PSL(2,\RR)$ be a non-cocompact Fuchsian group isomorphic to $\FF_3$, the free group on three generators. Consider the action of $\G_1$ on the boundary $S^1$ of the Poincar\'e disc. This action is $2$-filling and the algebra $\cA_1=C(S^1) \rtimes \G_1$ is p.i.s.u.n. with K-theory given by $K_0(\cA_1)=K_1(\cA_1)= \ZZ^4$, $[\Ind]=(1,0,0,0)$ \cite{ad2}. (The K-theory is independent of the embedding $\G_1 \subset \PSL(2,\RR)$.)

Let $\G_2$ be the $\tA_2$ group B.3 of \cite{cmsz}. This group is a lattice subgroup of 
$\PGL_3(\QQ_2)$ and acts naturally on the corresponding building of type $\tA_2$ and its boundary $\Om$. By Remark \ref{not2filling}, $2< \phi(\G,\Om) \le 6$.
By \cite{rs'}, the algebra $\cA_2=C(\Om) \rtimes \G_2$ is p.i.s.u.n. and satisfies the Universal Coefficient Theorem.  By \cite{rs''} the K-theory of $\cA_2$ is given by $K_0(\cA_2)=K_1(\cA_2)= \ZZ^4$, $[\Ind]=0$.

It follows from the classification theorem of \cite{k} that $\cA_1$, $\cA_2$ are stably isomorphic (but not isomorphic, since the classes $[\Ind]$ do not correspond).
}
\end{example}

\end{document}